\documentclass[12pt]{article}

\usepackage{amssymb,amsmath}
\newtheorem{theorem}{Theorem}[section]

\newtheorem{e-proposition}[theorem]{Proposition}

\newtheorem{e-definition}[theorem]{Definition\rm}

\def\og{\leavevmode\raise.3ex\hbox{$\scriptscriptstyle\langle\!\langle$~}}
\def\fg{\leavevmode\raise.3ex\hbox{~$\!\scriptscriptstyle\,\rangle\!\rangle$}}

\def\PP{\mathbb{P}}
\def\EE{\mathbb{E}}
\def\RR{\mathbb{R}}
\def\ZZ{\mathbb{Z}}
\def\eps{\varepsilon}
\def\a{{\bf a}}
\def\m{{\bf m}}
\def\v{\overrightarrow}
\newcommand\sy[1]{{#1}^{\#}}

\begin{document}
% place in the next line the header (rubrique) chosen for your article,
% if you know it (you can also have 2, format : Header1/Header2

\title{Bounds on the concentration function \\ 
in terms of Diophantine approximation}

\author{Omer Friedland and Sasha Sodin}

\maketitle
\footnotetext[1]{[omerfrie; sodinale]@post.tau.ac.il; address:
School of Mathematical Sciences, Tel Aviv University, Ramat Aviv,
Tel Aviv 69978, Israel}

\begin{abstract}

We demonstrate a simple analytic argument that may be used to bound the L\'evy
concentration function of a sum of independent random variables. The main application
is a version of a recent inequality due to Rudelson and Vershynin, and its
multidimensional generalisation.

\vskip 0.5\baselineskip

\noindent {\bf Des bornes pour la fonction de concentration en mati\`ere d'appro\-xi\-mation
Diophantienne.} Nous montrons un simple raisonnement analytique qui peut \^etre
utile pour borner la fonction de concentration d'une somme des variables
al\'eatoires ind\'ependantes. L'application principale est une version de l'in\'egalit\'e
r\'ecente de Rudelson et Vershynin, et sa g\'en\'eralisation au cadre multidimensionel.

\end{abstract}

\section{Introduction}

\noindent The {\em P.\ L\'evy concentration function} of a random variable $S$ is
defined as 
\[ \mathcal{Q}(S) = \sup_{x \in \RR} \, \PP \big\{ |S - x| \leq 1 \big\}~.\]
Since the work of L\'evy, Littlewood--Offord, Erd\H{o}s, Esseen, Kolmogorov and others,
numerous results in probability theory concern upper bounds on the concentration
function of the sum of independent random variables; a particularly powerful approach
was introduced in the 1970-s by Hal\'asz \cite{H}.

This note was motivated by the recent work of Rudelson and Vershynin \cite{RV}.
Let $X$ be a random variable; let $X_1, \cdots, X_n$ be independent copies
of $X$, and let $\a = (a_1, \cdots, a_n)$ be an $n$-tuple of real numbers.

In the Gaussian case $X \sim N(0, 1)$, we have: $\sum_{k=1}^n a_k X_k \sim N(0, |\a|^2)$
(where $|\cdot|$ stands for Euclidean norm), and consequently
\begin{equation}\label{g}
\mathcal{Q} \left( \sum_{k=1}^n a_k X_k \right)
    = \sqrt{\frac{2}{\pi|\a|}}  \, (1 + o(1))~, \quad |\a| \to \infty~.
\end{equation}
On the other hand, if $X$ has atoms, the left-hand side of (\ref{g}) does not tend
to $0$ as $|\a| \to \infty$. Therefore one may ask, for which $\a \in \RR^n$ is it true that
\begin{equation}\label{leqdel}
\mathcal{Q} \left( \sum_{k=1}^n a_k X_k \right) \leq C/|\a|~?
\end{equation}

Rudelson and Vershynin gave a bound in terms of Diophantine approximation of the
vector $\a$. Their approach makes use of a deep measure-theoretic lemma from \cite{H}.
Our goal is to show a simpler analytic method that may be of use in such problems.
The following theorem is a (slightly improved) version of \cite[Theorem 1.3]{RV}.
\vspace{1mm}
\begin{theorem}\label{thm1} Let ${X}_1, \cdots, {X}_n$ be independent copies of a random
variable ${X}$ such that $\mathcal{Q}(X) \leq 1 - p$,
and let $\a = (a_1, \cdots, a_n) \in \RR^n$. If, for some $0 < D < 1$ and $\alpha > 0$,
\begin{equation}\label{ass1}
|\eta \a - \m| \geq \alpha
    \quad \text{for} \quad
    \m \in \ZZ^n, \, \eta \in \big[1/(2\| a \|_\infty) , \,  D\big]~,
\end{equation}
then
\begin{equation}\label{conc1}
\mathcal{Q}(\sum_{k=1}^n X_k a_k)
    \leq C \left\{ \exp(-c p^2 \alpha^2) + \frac{1}{p \, D} \, \frac{1}{|\a|} \right\}~.
\end{equation}
\end{theorem}
Here and further $C,c, C', c_1,\cdots>0$ denote numerical constants. \hfill \vspace{2mm}\\
We also extend this result to the multidimensional case. The concentration function of an
$\RR^d$-valued random vector $\v{S}$ is defined as
\[ \mathcal{Q}(\v{S}) = \sup_{x \in \RR^d} \, \PP \big\{ |\v{S} - \v{x}| \leq 1 \big\}~.\]
\vspace{1mm}
\begin{theorem}\label{thm} Let ${X}_1, \cdots, {X}_n$ be independent copies of a
random variable ${X}$ such that $\mathcal{Q}(X) \leq 1-p$, and let
$\v{a}_1, \cdots, \v{a}_n \in \RR^d$ be such that, for some $0 < D < d$ and $\alpha > 0$, %%%
\begin{equation}\label{assum}\begin{split}
\sum_{k=1}^n (\v\eta \cdot \v{a}_k - m_k)^2 \geq \alpha^2
\,\,\, &\text{for} \,\,\, m_1,\cdots,m_n \in \ZZ , \, \v\eta \in \RR^d \\
\,\,\, &\text{such that} \,\,\,
\max_k |\v\eta \cdot \v{a}_k| \geq 1/2, \, |\v\eta| \leq D.
\end{split}\end{equation}
Then
\begin{equation}\label{conc}\begin{split}
\mathcal{Q}(\sum_{k=1}^n X_k \v{a}_k)
    &\leq C^d \Big\{  \exp(-c p^2 \alpha^2) \\
    &\qquad+ \left( \frac{\sqrt{d}}{pD} \right)^d
        \det{}^{-1/2} \left[ \sum_{k=1}^n \v{a}_k \otimes \v{a}_k \right] \Big\}~,
\end{split}\end{equation}
where $C,c>0$ are numerical constants.
\end{theorem}
\vspace{1mm}
Of course, Theorem~\ref{thm1} follows formally from Theorem~\ref{thm}. For simplicity
of exposition we will prove Theorem~\ref{thm1} and indicate the adjustments that are necessary
for $d > 1$.

\section{Proof of Theorem~\ref{thm1}}

\noindent {\bf Step 1}: By Chebyshev's inequality and the identity
\begin{equation*}
\exp(-y^2) = \int_{-\infty}^{+\infty}
    \exp \left\{ 2iy\eta - \eta^2 \right\} \frac{d\eta}{\sqrt\pi}
\end{equation*}
it follows that
\[\begin{split}
\PP \left\{ \left| \sum_{k=1}^n X_k a_k - x \right| \leq 1 \right\}
    &\leq e \, \EE \exp \left\{ - \left[ \sum_{k=1}^n X_k a_k - x\right]^2 \right\} \\
    &= e \, \EE \int_{-\infty}^{+\infty} \exp \left\{ 2 i \left[ \sum_{k=1}^n X_k a_k - x \right] \eta - \eta^2 \right\}
        \frac{d\eta}{\sqrt{\pi}}~.
\end{split}\]
Now we can swap the expectation with the integral and take absolute value:
\begin{equation*}\begin{split}
\PP \left\{ \left| \sum_{k=1}^n X_k a_k - x \right| \leq 1 \right\}
    &\leq e \int_{-\infty}^{+\infty} \prod_{k=1}^n \phi(2 a_k \eta)
        \exp \left\{ - 2ix\eta - \eta^2 \right\} \frac{d\eta}{\sqrt\pi} \\
    &\leq e \int_{-\infty}^{+\infty} \prod_{k=1}^n |\phi(2 a_k \eta)|
        \exp \left\{ - \eta^2 \right\} \frac{d\eta}{\sqrt\pi}~,
\end{split}\end{equation*}
where $\phi(\eta) = \EE \exp(i \eta X)$ is the characteristic function
of every one of the $X_k$. Therefore
\begin{equation}\label{step1}
\mathcal{Q}(\sum_{k=1}^n X_k a_k)
    \leq e \int_{-\infty}^{+\infty} \prod_{k=1}^n |\phi(2 a_k \eta)|
        \exp \left\{ - \eta^2 \right\} \frac{d\eta}{\sqrt\pi}~.
\end{equation}

\vspace{2mm}
\noindent {\bf Step 2} (this step is analogous to \cite[\S 3]{H} and \cite[4.2]{RV}):
First,
\[ |\phi(\eta)| \leq \exp \left( - \frac{1}{2} (1 - |\phi(\eta)|^2) \right)~.\]
Let $X'$ be an independent copy of $X$, $\sy{X} = X-X'$. Observe that
\[ q = \PP \left\{ |\sy{X}|\geq 2 \right\} \geq p^2/2 \]
and
\[\begin{split}
|\phi(\eta)|^2
   &= \EE \, \exp \big(i \eta \sy{X} \big) \\
   &= \EE  \cos\big(\eta \sy{X} \big)
   \leq (1 - q) + q \EE \left\{ \cos\big(\eta \sy{X} \big) \, \big| \, |\sy{X}| \geq 2 \right\}~; 
\end{split}\]
therefore
\begin{multline*}
\int_{-\infty}^{+\infty} \prod_{k=1}^n |\phi(2 a_k \eta)|
        \exp \left\{ - \eta^2 \right\} \frac{d\eta}{\sqrt\pi} \\
    \leq \int_{-\infty}^{+\infty}
        \exp \left\{ - \frac{q}{2} \,
            \EE \left[ \sum_{k=1}^n \left( 1 - \cos \big( 2 a_k \eta \sy{X} \big)\right)
                \, \Big| \, |\sy{X}| \geq 2 \right] - \eta^2 \right\}
        \frac{d\eta}{\sqrt\pi} \\
    \leq \EE  \left[ \int_{-\infty}^{+\infty}
        \exp \left\{ - \frac{q}{2} \,
            \sum_{k=1}^n \left( 1 - \cos \big( 2 a_k \eta \sy{X} \big)\right) - \eta^2 \right\}
        \frac{d\eta}{\sqrt\pi} \,\, \Big| \,\, |\sy{X}| \geq 2 \right]~.
\end{multline*}

Replace the conditional expectation with supremum over the possible values $z \geq 2$ of
$|\sy{X}|$ and recall that
\[ 1 - \cos \theta \geq c_1 \min_{m \in \ZZ} |\theta - 2 \pi m|^2~; \]
then
\begin{equation}\label{step2}\begin{split}
&\int_{-\infty}^{+\infty} \prod_{k=1}^n |\phi(2 \eta a_k )|
    \exp \left\{ - \eta^2 \right\} \frac{d\eta}{\sqrt\pi} \\
    &\qquad\leq \sup_{z \geq 2}
        \int_{-\infty}^{+\infty} \exp \left\{ - \frac{q}{2} \,
            \sum_{k=1}^n \left( 1 - \cos \big( 2 z \eta a_k  \big)\right) - \eta^2 \right\}
        \frac{d\eta}{\sqrt\pi} \\
    &\qquad\leq \sup_{z \geq 2}
        \int_{-\infty}^{+\infty} \exp \left\{ - c_2 p^2 \,
            \sum_{k=1}^n \min_{m_k} |z \eta a_k - \pi m_k|^2 - \eta^2 \right\}
        \frac{d\eta}{\sqrt\pi} \\
    &\qquad= \sup_{z \geq 2/\pi}
        \int_{-\infty}^{+\infty} \exp \left\{ - c_3 p^2 \,
            \sum_{k=1}^n \min_{m_k} |\eta a_k - m_k|^2 - (\eta/z)^2 \right\}
        \frac{d\eta}{z\sqrt\pi}~.
\end{split}\end{equation}

\vspace{2mm}
\noindent {\bf Step 3:} Denote
\[ A = \left\{ \eta \in \RR \, \Big| \,
    \forall \, \m \in \ZZ^n \, : \,
    |\eta \a  - \m| \geq \alpha/2 \right\} \, , \quad B = \RR \setminus A~. \]
Then the last integral in (\ref{step2}) can be split into
\begin{equation}\label{triv}
\int_{-\infty}^{+\infty} = \int_A + \int_B~,
\end{equation}
and
\begin{equation}\label{estA}
\int_{A} \leq \exp \left( - c_3 p^2 \alpha^2\right)~.
\end{equation}
On the other hand, if $\eta', \eta'' \in B$, then
$|\eta' \a  - \pi \m'|, |\eta'' \a  - \pi \m''|< \alpha/2$ for some $\m', \m'' \in \ZZ^n$,
and hence
\[ \left| (\eta' - \eta'') \a - (\m' - \m'') \right| < \alpha~.\]
Therefore by (\ref{ass1}) either $|\eta' - \eta''| < 1 / (2 \| a \|_\infty)$
or $|\eta' - \eta''| > D$. In other words, $B \subset \bigcup_j B_j$, where $B_j$ are
intervals of length $\leq 1/\|a\|_\infty$ such that any two points belonging to
different $B_j$ are at least $D$-apart.

\vspace{2mm}
\noindent {\bf Step 4:} For every $j$ there exists $\eta_j \in B_j$ such that
\[ \int_{B_j}
    = \exp(-\eta_j^2/z^2) \int_{B_j} \exp \left\{ - c_3 p^2 \,
        \sum_{k=1}^n \min_{m_k} |\eta a_k - m_k|^2 \right\}
        \frac{d\eta}{z\sqrt\pi}~.\]
By H\"older's inequality
\begin{equation}\label{hoeld}
\int_{B_j} \leq e^{-\eta_j^2/z^2} \prod_{k=1}^n \left\{
    \int_{B_j} \exp \left\{ - \frac{c_3 p^2 |\a|^2}{a_k^2} \,
        \min_{m_k} |\eta a_k - m_k|^2 \right\}
        \frac{d\eta}{z\sqrt\pi}\right\}^\frac{a_k^2}{|\a|^2}~.
\end{equation}
The length of the interval $a_k B_j$ is $\leq 1$; hence $m_k$ (which is the closest integer
to $\eta a_k$) can obtain at most 2 values while $\eta \in B_j$. Therefore every one of
the integrals on the right-hand side of (\ref{hoeld}) is bounded by
\[ 2 \int_{-\infty}^{+\infty}
    \exp \left\{ - c_3 p^2 |\a|^2 \eta^2 \right\}
        \frac{d\eta}{z\sqrt\pi} = \frac{C_1}{z p |\a|}~,   \]
and therefore
\[ \int_{B} \leq \sum_j \int_{B_j} \leq \frac{C_1}{z p |\a|} \, \sum_j \exp(-\eta_j^2/z^2)~.\]
Now, $B_j$ (and hence $\eta_j$) are $D$-separated; therefore
\[\begin{split}
\sum_j \exp(-\eta_j^2/z^2)
    &\leq 2 \sum_{j=0}^\infty  \exp \left(- (Dj/z)^2 \right) \\
    &\leq 2 \left\{1 + \int_0^{+\infty} \exp \left( - (D \eta/z)^2 d\eta \right) \right\} \\
    &\leq 2 (1 + C_2z/D) \leq C_3z/D~.
\end{split}\]
Hence finally
\[ \int_B \leq \frac{C_4}{pD|\a|}~;\]
combining  this with (\ref{step1}--\ref{estA}) we deduce (\ref{conc1}).

\section{Remarks}

\begin{enumerate}
\item The results can be also used to estimate the formally more general form of the
L\'evy concentration function:
\[ \mathcal{L}(\sum_{k=1}^n X_k \v{a}_k; \eps)
    = \sup_{\v{x} \in \RR^d} \, \PP \big\{ |\sum_{k=1}^n X_k \v{a}_k - \v{x}| \leq \eps \big\}~. \]
Indeed, $\mathcal{L}(\sum_{k=1}^n X_k \v{a}_k; \eps) = \mathcal{Q}(\sum_{k=1}^n X_k \v{a}_k/\eps)$, so one can
just apply the theorems to $\v{a}_k' = \v{a}_k / \eps$. \vspace{1mm}
\item By similar reasoning, the assumption $\mathcal{Q}(X) \leq 1 - p$ can be replaced with
$\mathcal{Q}(K X) \leq 1 - p$ (for an arbitrary $K>0$); this will only influence the values of
the constants in (\ref{conc1}), (\ref{conc}). \vspace{1mm}
\item The proof of Theorem~\ref{thm} is parallel to that of Theorem~\ref{thm1}. The main
difference appears in Step~4, where instead of H\"older's inequality one should use
the Brascamp--Lieb--Luttinger rearrangement inequality \cite{BLL}. (Note that a different rearrangement
inequality was applied to a similar problem by Howard \cite{Ho}).
\end{enumerate}

\section*{Acknowledgements} We are grateful to our supervisor Vitali Milman for his support
and for encouraging to write this note. We thank Mark Rudelson and Roman Vershynin for
stimulating discussions, and in particular for suggesting the current formulation of
Theorem~\ref{thm} with improved dependence on the dimension $d$, and for spotting several
blunders.

\end{document}